# Upper bound for the second Hankel determinant of certain subclass of analytic and bi-univalent functions

Nizami Mustafa

Department of Mathematics, Faculty of Science and Letters, Kafkas University, Kars, Turkey

**Abstract**

In this paper, we consider a general subclass of analytic and bi-univalent functions in the open unit disk in the complex plane. Making use of the Chebyshev polynomials, we obtain upper bound estimate for the second Hankel determinant for this function class.

## 1. Introduction and preliminaries

Let $A$ denote the class of all complex valued functions $f(z)$ given by

$$f(z) = z + a_2 z^2 + a_3 z^3 + \cdots = z + \sum_{n=2}^{\infty} a_n z^n, \quad (1.1)$$

which are analytic in the open unit disk $U = \{z \in \mathbb{C} : |z| < 1\}$ in the complex plane. Furthermore, let $S$ be the class of all functions in $A$ which are univalent in $U$. Some of the important and well-investigated subclasses of $S$ include the classes $S^*(\alpha)$ and $C(\alpha)$ given below (see also [6, 8, 19]) such that $S^*(\alpha)$ is the class of starlike functions of order $\alpha$ and $C(\alpha)$ is the class of convex functions of order $\alpha$ ($\alpha \in [0,1)$)

$$S^*(\alpha) = \left\{ f \in S : \operatorname{Re}\left[\frac{zf'(z)}{f(z)}\right] > \alpha, z \in U \right\}, \; \alpha \in [0,1)$$

and

$$C(\alpha) = \left\{ f \in S : \operatorname{Re}\left[1 + \frac{zf''(z)}{f'(z)}\right] > \alpha, z \in U \right\}, \; \alpha \in [0,1).$$

For $f \in S$, we consider the following class of analytic functions

$$Q(\alpha, \beta) = \left\{ f \in S : \operatorname{Re}\left[(1-\beta)\frac{f(z)}{z} + \beta f'(z)\right] > \alpha, z \in U \right\}, \; \alpha \in [0,1), \beta \in [0,1].$$

In 1995, Ding et al. [4] have introduced and investigated the class $Q(\alpha,\beta)$. It is well-known that (see [6]) every function $f \in S$ has an inverse $f^{-1}$, defined by

$$f^{-1}(f(z)) = z,\ z \in U$$

and

$$f(f^{-1}(w)) = w,\ w \in D_{r_0} = \{w \in \mathbb{C}:\ |w| < r_0(f)\},\ r_0(f) \geq 1/4,$$

where

$$f^{-1}(w) = w - a_2 w^2 + (2a_2^2 - a_3)w^3 - (5a_2^3 - 5a_2 a_3 + a_4)w^4 + \cdots,\ w \in D_{r_0}. \tag{1.2}$$

A function $f \in A$ is called bi-univalent in $U$ if both $f$ and $f^{-1}$ are univalent in the definition sets. Let $\Sigma$ denote the class of bi-univalent functions in $U$ given by (1.1). For a short history and examples of functions in the class $\Sigma$, see [18].

Firstly, Lewin [16] introduced the class of bi-univalent functions, obtaining the estimate $|a_2| \leq 1.51$. Subsequently, Brannan and Clunie [2] developed the result of Lewin to $|a_2| \leq \sqrt{2}$ for $f \in \Sigma$. Accordingly, Netanyahu [14] showed that $|a_2| \leq \frac{4}{3}$. Earlier, Brannan and Taha [1] introduced certain subclasses of bi-univalent function class $\Sigma$, namely bi-starlike function of order $\alpha$ denoted $S_\Sigma^*(\alpha)$ and bi-convex function of order $\alpha$ denoted $C_\Sigma(\alpha)$ corresponding to the function classes $S^*(\alpha)$ and $C(\alpha)$, respectively. For each of the function classes $S_\Sigma^*(\alpha)$ and $C_\Sigma(\alpha)$, non-sharp estimates on the first two Taylor-Maclaurin coefficients were found in [1, 21]. Many researchers (see [20, 22, 23]) have introduced and investigated several interesting subclasses of bi-univalent function class $\Sigma$ and they have found non-sharp estimates on the first two Taylor-Maclaurin coefficients. However, the coefficient problem for each of the Taylor-Maclaurin coefficients $|a_n|$, $n = 3, 4, \ldots$ is still an open problem (see, for example, [14, 16]).

An analytic function $f$ is bi-starlike of Ma-Minda type or bi-convex of Ma-Minda type if both $f$ and $f^{-1}$ are, respectively, Ma-Minda starlike and convex. These classes are denoted, respectively, by $S_\Sigma^*(\phi)$ and $C_\Sigma(\phi)$. In the sequel, it is assumed that the function $\phi$ is an analytic function with positive real part in $U$, satisfying $\phi(0) = 1$, $\phi'(0) > 0$ and $\phi(U)$ is symmetric with respect to the real axis. Such a function has a series expansion of the following form:

$$\phi(z) = 1 + b_1 z + b_2 z^2 + b_3 z^3 + \cdots = 1 + \sum_{n=1}^{\infty} b_n z^n,\ b_1 > 0. \tag{1.3}$$

An analytic function $f$ is subordinate to an analytic function $\phi$, written $f(z) \prec \phi(z)$, provided that there is an analytic function (that is, Schwarz function) $\omega$ defined on $U$ with $\omega(0) = 0$ and $|\omega(z)| < 1$ satisfying $f(z) = \phi(\omega(z))$. Ma and Minda [12] unified various subclasses of starlike and convex functions for which either of the quantity $\frac{zf'(z)}{f(z)}$ or $1+\frac{zf''(z)}{f'(z)}$ is subordinate to a more general function. For this purpose, they considered an analytic function $\phi$ with positive real part in $U$, $\phi(0)=1$, $\phi'(0)>0$ and $\phi$ maps $U$ onto a region starlike with respect to 1 and symmetric with respect to the real axis. The class of Ma-Minda starlike and Ma-Minda convex functions consists of functions $f \in A$ satisfying the subordination $\frac{zf'(z)}{f(z)} \prec \phi(z)$ and $1+\frac{zf''(z)}{f'(z)} \prec \phi(z)$, respectively.

In 1976, Noonan and Thomas [15] defined the $q$th Hankel determinant of $f$ for $q \in \mathbb{N}$ by

$$H_q(n) = \begin{vmatrix} a_n & \cdots & a_{n+q-1} \\ \cdot & \cdots & \cdot \\ a_{n+q-1} & \cdots & a_{n+2q-2} \end{vmatrix}.$$

For $q=2$ and $n=1$, Fekete and Szegö [7] considered the Hankel determinant of $f$ as $H_2(1) = \begin{vmatrix} a_1 & a_2 \\ a_2 & a_3 \end{vmatrix} = a_1 a_3 - a_2^2$. They made an earlier study for the estimates of $|a_3 - \mu a_2^2|$ when $a_1 = 1$ with real $\mu \in \mathbb{R}$. The well-known result due to them states that if $f \in A$, then

$$|a_3 - \mu a_2^2| \leq \begin{cases} 3-4\mu & \text{if } \mu \in (-\infty, 0], \\ 1+2\exp\left(\frac{-2\mu}{1-\mu}\right) & \text{if } \mu \in [0,1], \\ 4\mu-3 & \text{if } \mu \in [1,+\infty). \end{cases}$$

Furthermore, Hummel [9, 10] obtained sharp estimates for $|a_3 - \mu a_2^2|$ when $f$ is a convex function and also Keogh and Merkes [11] obtained sharp estimates for $|a_3 - \mu a_2^2|$ when $f$ is a close-to-convex function, starlike and convex function in $U$.

Recently, the upper bounds of $|H_2(2)| = |a_2 a_4 - a_3^2|$ for the classes $S_\Sigma^*(\alpha)$ and $C_\Sigma(\alpha)$ were obtained by Deniz et al. [3]. Very soon, Orhan et al. [17] reviwed the study of bounds for the second Hankel determinant of the subclass $M_\Sigma^\alpha(\beta)$ of bi-univalent functions.

Chebyshev polynomials, which is used by us in this paper, play a considerable act in numerical analysis and mathematical physics. It is well-known that the Chebyshev polinomials are four kinds. The most of research articles related to specific orthogonal polynomials of Chebyshev family, contain essentially results of Chebyshev polynomials of first and second kinds $T_n(x)$ and $U_n(x)$, and their numerous uses in different applications (see [5, 13]).

The well-known kinds of the Chebyshev polynomials are the first and second kinds. In the case of real variable $x$ on $(-1,1)$, the first and second kinds of the Chebyshev polynomials are defined by

$$T_n(x) = \cos(n \arccos x),$$

$$U_n(x) = \frac{\sin[(n+1)\arccos x]}{\sin(\arccos x)} = \frac{\sin[(n+1)\arccos x]}{\sqrt{1-x^2}}.$$

We consider the function

$$G(t,z) = \frac{1}{1-2tz+z^2}, \ t \in \left(\frac{1}{2}, 1\right), z \in U.$$

It is well-known that if $t = \cos\alpha, \ \alpha \in \left(0, \frac{\pi}{3}\right)$, then

$$G(t,z) = 1 + \sum_{n=1}^{\infty} \frac{\sin[(n+1)\alpha]}{\sin\alpha} z^n =$$
$$1 + 2\cos\alpha z + (3\cos^2\alpha - \sin^2\alpha)z^2 + (8\cos^3\alpha - 4\cos\alpha)z^3 + \cdots, \ z \in U.$$

That is,

$$G(t,z) = 1 + U_1(t)z + U_2(t)z^2 + U_3(t)z^3 + \cdots, \ t \in \left(\frac{1}{2}, 1\right), z \in U, \qquad (1.4)$$

where $U_n(t) = \frac{\sin[(n+1)\arccos t]}{\sqrt{1-t^2}}$, $n \in \mathbb{N}$ are the second kind Chebyshev polynomials. From the definition of the second kind Chebyshev polynomials, we easily obtain that $U_1(t) = 2t$. Also, it is well-known that

$$U_{n+1}(t) = 2tU_n(t) - U_{n-2}(t)$$

for all $n \in \mathbb{N}$. From here, we can easily obtain

$$U_2(t) = 4t^2 - 1, \ U_3(t) = 8t^3 - 4t. \qquad (1.5)$$

Inspired by the aforementioned works, making use of the Chebyshev polynomials, we define a subclass of bi-univalent functions $\Sigma$ as follows.

**Definition 1.1.** *A function $f \in \Sigma$ given by (1.1) is said to be in the class $Q_\Sigma(G, \beta, t)$, $\beta \in [0,1], t \in \left(\frac{1}{2}, 1\right)$, where $G$ is an analytic function given by (1.4), if the following conditions are satisfied*

$$(1-\beta)\frac{f(z)}{z}+\beta f'(z) \prec G(t,z),\ z \in U \qquad (1.6)$$

and

$$(1-\beta)\frac{g(w)}{w}+\beta g'(w) \prec G(t,w),\ w \in D_{r_0} \qquad (1.7)$$

where $g = f^{-1}$.

**Remark 1.1.** *Taking* $\beta = 1$, *we have* $Q_\Sigma(G,1,t) = \Re_\Sigma(G,t),\ t \in \left(\frac{1}{2},1\right)$; *that is,*

$$f \in \Re_\Sigma(G,t) \Leftrightarrow f'(z) \prec G(t,z),\ z \in U\ and\ g'(w) \prec G(t,w),\ w \in D_{r_0},$$

where $g = f^{-1}$.

**Remark 1.2.** *Taking* $\beta = 0$, *we have* $Q_\Sigma(G,0,t) = N_\Sigma(G,t),\ t \in \left(\frac{1}{2},1\right)$; *that is,*

$$f \in N_\Sigma(G,t) \Leftrightarrow \frac{f(z)}{z} \prec G(t,z),\ z \in U\ and\ \frac{g(w)}{w} \prec G(t,w),\ w \in D_{r_0},$$

where $g = f^{-1}$.

The object of this paper is to determine the second Hankel determinant for the function class $Q_\Sigma(G,\beta,t)$ and its special classes, and is to give upper bound estimate for $|H_2(2)|$. In the study, the bound estimates for the initial coefficients of the functions belonging to this class are also obtained. For the functions belonging to class $Q_\Sigma(G,\beta,t)$, Fekete-Szöge inequality is also obtained.

In order to prove our main results, we shall need the following lemma.

**Lemma 1.1.** ([4]) *Let* P *be the class of all analytic functions* $p(z)$ *of the form*

$$p(z) = 1 + p_1 z + p_2 z^2 + \cdots = 1 + \sum_{n=1}^{\infty} p_n z^n \qquad (1.8)$$

*satisfying* $\operatorname{Re}(p(z)) > 0,\ z \in U$ *and* $p(0) = 1$. *Then,* $|p_n| \leq 2$, *for every* $n = 1,2,3,\ldots$ . *This inequality is sharp for each* $n$.
*Moreover,*

$$2p_2 = p_1^2 + (4-p_1^2)x,$$
$$4p_3 = p_1^3 + 2(4-p_1^2)p_1 x - (4-p_1^2)p_1 x^2 + 2(4-p_1^2)(1-|x|^2)z,$$

*for some* $x$, $z$ *with* $|x| \leq 1$, $|z| \leq 1$.

## 2. Upper bound for the second Hankel determinant of the class $Q_\Sigma(G, \beta, t)$

In this section, we prove the following theorem on upper bound of the second Hankel determinant of the function class $Q_\Sigma(G, \beta, t)$.

**Theorem 2.1.** *Let the function* $f(z)$ *given by (1.1) be in the class* $Q_\Sigma(G, \beta, t)$, $\beta \in [0,1], t \in \left(\frac{1}{2}, 1\right)$, *where the function* $G$ *is an analytic function given by (1.4). Then,*

$$|a_2 a_4 - a_3^2| \leq \begin{cases} H(t, 2-), & \text{if } \Delta(\beta,t) \geq 0 \text{ and } c(\beta,t) \geq 0, \\ \max\left\{\dfrac{4t^2}{(1+2\beta)^2}, H(t, 2-)\right\}, & \text{if } \Delta(\beta,t) > 0 \text{ and } c(\beta,t) < 0, \\ \dfrac{4t^2}{(1+2\beta)^2}, & \text{if } \Delta(\beta,t) \leq 0 \text{ and } c(\beta,t) \leq 0, \\ \max\{H(t, \tau_0), H(t, 2-)\}, & \text{if } \Delta(\beta,t) < 0 \text{ and } c(\beta,t) > 0, \end{cases}$$

*where*

$$H(t, 2-) = \frac{8t^2 \left|(2t^2-1)(1+\beta)^3 - 2t^2(1+3\beta)\right|}{(1+\beta)^4(1+3\beta)},$$

$$H(t, \tau_0) = \frac{4t^2}{(1+2\beta)^2} - \frac{c^2(\beta,t)}{4(1+\beta)^4(1+2\beta)^2(1+3\beta)}, \quad p_0 = \sqrt{\frac{-2c(\beta,t)}{\Delta(\beta,t)}},$$

$$\Delta(\beta,t) = 16t^2 \left|(2t^2-1)(1+\beta)^3 - 2t^2(1+3\beta)\right|(1+2\beta)^2 -$$
$$8t\left[t^2 - (4t^2-1)(1+\beta)(1+2\beta)\right](1+\beta)^2(1+2\beta)(1+3\beta) - 8t^2(1+\beta)^3\beta^2,$$

$$c(\beta,t) = 8t\left[(5t^2-1)(1+3\beta) + 2(4t^2-1)\beta^2\right](1+\beta)^2(1+2\beta) -$$
$$8t^2\left[(1+\beta)^2 + (2+\beta)\beta\right](1+\beta)^3.$$

**Proof.** Let $f \in Q_\Sigma(G, \beta, t)$, $\beta \in [0,1], t \in \left(\frac{1}{2}, 1\right)$ and $g = f^{-1}$. Then, according to Definition 1.1, there are analytic functions $\omega: U \to U$, $\varpi: D_{r_0} \to D_{r_0}$ with $\omega(0) = 0 = \varpi(0)$, $|\omega(z)| < 1$, $|\varpi(w)| < 1$ satisfying the following conditions

$$(1-\beta)\frac{f(z)}{z} + \beta f'(z) = G(t, \omega(z)), \quad z \in U \tag{2.1}$$

and

$$(1-\beta)\frac{g(w)}{w} + \beta g'(w) = G(t, \varpi(w)), \quad w \in D_{r_0}. \tag{2.2}$$

Let also the functions $p, q \in P$ be define as follows

$$p(z) := \frac{1+\omega(z)}{1-\omega(z)} = 1 + p_1 z + p_2 z^2 + \cdots = 1 + \sum_{n=1}^{\infty} p_n z^n$$

and

$$q(w) := \frac{1+\varpi(w)}{1-\varpi(w)} = 1 + q_1 w + q_2 w^2 + \cdots = 1 + \sum_{n=1}^{\infty} q_n w^n.$$

It follows that

$$\omega(z) := \frac{p(z)-1}{p(z)+1} = \frac{1}{2}\left[p_1 z + \left(p_2 - \frac{p_1^2}{2}\right)z^2 + \left(p_3 - p_1 p_2 + \frac{p_1^3}{4}\right)z^3 + \cdots\right] \quad (2.3)$$

and

$$\varpi(w) := \frac{q(w)-1}{q(w)+1} = \frac{1}{2}\left[q_1 w + \left(q_2 - \frac{q_1^2}{2}\right)w^2 + \left(q_3 - q_1 q_2 + \frac{q_1^3}{4}\right)w^3 + \cdots\right]. \quad (2.4)$$

From (2.3) and (2.4), considering (1.4), we can easily show that

$$G(t, \omega(z)) = 1 + \frac{U_1(t)}{2} p_1 z + \left[\frac{U_1(t)}{2}\left(p_2 - \frac{p_1^2}{2}\right) + \frac{U_2(t)}{4} p_1^2\right] z^2 +$$
$$\left[\frac{U_1(t)}{2}\left(p_3 - p_1 p_2 + \frac{p_1^3}{4}\right) + \frac{U_2(t)}{2} p_1 \left(p_2 - \frac{p_1^2}{2}\right) + \frac{U_3(t)}{8} p_1^3\right] z^3 + \cdots \quad (2.5)$$

and

$$G(t, \varpi(w)) = 1 + \frac{U_1(t)}{2} q_1 w + \left[\frac{U_1(t)}{2}\left(q_2 - \frac{q_1^2}{2}\right) + \frac{U_2(t)}{4} q_1^2\right] w^2 +$$
$$\left[\frac{U_1(t)}{2}\left(q_3 - q_1 q_2 + \frac{q_1^3}{4}\right) + \frac{U_2(t)}{2} q_1 \left(q_2 - \frac{q_1^2}{2}\right) + \frac{U_3(t)}{8} q_1^3\right] w^3 + \cdots. \quad (2.6)$$

From (2.1), (2.5) and (2.2), (2.6), we can easily obtain that

$$(1+\beta)a_2 = \frac{U_1(t)}{2} p_1, \quad (2.7)$$

$$(1+2\beta)a_3 = \frac{U_1(t)}{2}\left(p_2 - \frac{p_1^2}{2}\right) + \frac{U_2(t)}{4} p_1^2, \quad (2.8)$$

$$(1+3\beta)a_4 = \frac{U_1(t)}{2}\left(p_3 - p_1 p_2 + \frac{p_1^3}{4}\right) + \frac{U_2(t)}{2} p_1 \left(p_2 - \frac{p_1^2}{2}\right) + \frac{U_3(t)}{8} p_1^3 \quad (2.9)$$

and

$$-(1+\beta)a_2 = \frac{U_1(t)}{2}q_1, \tag{2.10}$$

$$(1+2\beta)(2a_2^2 - a_3) = \frac{U_1(t)}{2}\left(q_2 - \frac{q_1^2}{2}\right) + \frac{U_2(t)}{4}q_1^2, \tag{2.11}$$

$$-(1+3\beta)(5a_2^3 - 5a_2a_3 + a_4) = \frac{U_1(t)}{2}\left(q_3 - q_1q_2 + \frac{q_1^3}{4}\right) + \frac{U_2(t)}{2}q_1\left(q_2 - \frac{q_1^2}{2}\right) + \frac{U_3(t)}{8}q_1^3. \tag{2.12}$$

From (2.7) and (2.10), we obtain that

$$\frac{U_1(t)}{2(1+\beta)}p_1 = a_2 = -\frac{U_1(t)}{2(1+\beta)}q_1. \tag{2.13}$$

Subtracting (2.11) from (2.8) and considering (2.13), we can easily obtain that

$$a_3 = a_2^2 + \frac{U_1(t)}{4(1+2\beta)}(p_2 - q_2) = \frac{U_1^2(t)}{4(1+\beta)^2} + \frac{U_1(t)}{4(1+2\beta)}(p_2 - q_2). \tag{2.14}$$

On the other hand, subtracting (2.12) from (2.9) and considering (2.13) and (2.14), we get

$$a_4 = \frac{5U_1^2(t)p_1(p_2 - q_2)}{6(1+\beta)(1+2\beta)} + \frac{U_1(t)(p_3 - q_3)}{4(1+3\beta)} + \frac{U_2(t) - U_1(t)}{4(1+3\beta)}(p_2 + q_2) + \frac{U_1(t) - 2U_2(t) + U_3(t)}{8(1+3\beta)}p_1^3. \tag{2.15}$$

Thus, from (2.13), (2.14) and (2.15), we can easily establish that

$$a_2a_4 - a_3^2 = \frac{U_1^3(t)p_1^2(p_2 - q_2)}{32(1+\beta)^2(1+2\beta)} + \frac{U_1^2(t)p_1(p_3 - q_3)}{8(1+\beta)(1+3\beta)} + \frac{[U_2(t) - U_1(t)]U_1(t)}{8(1+\beta)(1+3\beta)}p_1^2(p_2 + q_2) - \frac{U_1^2(t)(p_2 - q_2)^2}{16(1+2\beta)^2} + \frac{U_1(t)p_1^4}{16(1+\beta)^4(1+3\beta)}\left\{[U_1(t) - 2U_2(t) + U_3(t)](1+\beta)^3 - U_1^3(t)(1+3\beta)\right\}. \tag{2.16}$$

According to Lemma 1.1, we have

$$2p_2 = p_1^2 + (4 - p_1^2)x \text{ and } 2q_2 = q_1^2 + (4 - q_1^2)y \tag{2.17}$$

and

$$4p_3 = p_1^3 + 2(4-p_1^2)p_1 x - (4-p_1^2)p_1 x^2 + 2(4-p_1^2)(1-|x|^2)z$$
$$4q_3 = q_1^3 + 2(4-q_1^2)q_1 y - (4-q_1^2)q_1 y^2 + 2(4-q_1^2)(1-|y|^2)w, \quad (2.18)$$

for some $x, y, z, w$ with $|x| \le 1$, $|y| \le 1$, $|z| \le 1$, $|w| \le 1$.

Since (see (2.13)) $p_1 = -q_1$, from (2.17) and (2.18), we get

$$p_2 - q_2 = \frac{4-p_1^2}{2}(x-y), \quad p_2 + q_2 = p_1^2 + \frac{4-p_1^2}{2}(x+y) \quad (2.19)$$

and

$$p_3 - q_3 = \frac{p_1^3}{2} + \frac{(4-p_1^2)p_1}{2}(x+y) - \frac{(4-p_1^2)p_1}{4}(x^2+y^2) + \frac{4-p_1^2}{2}\left[(1-|x|^2)z - (1-|y|^2)w\right]. \quad (2.20)$$

According to Lemma 1.1, we may assume without any restriction that $\tau \in [0, 2]$, where $\tau = |p_1|$.

Thus, substituting the expressions (2.19) and (2.20) in (2.16) and using triangle inequality, letting $|x| = \xi$, $|y| = \eta$, we can easily obtain that

$$|a_2 a_4 - a_3^2| \le c_1(t,\tau)(\xi+\eta)^2 + c_2(t,\tau)(\xi^2+\eta^2) + c_3(t,\tau)(\xi+\eta) + c_4(t,\tau) := F(\xi,\eta), \quad (2.21)$$

where

$$c_1(t,\tau) = \frac{U_1^2(t)(4-\tau^2)^2}{64(1+2\beta)^2} \ge 0, \quad c_2(t,\tau) = \frac{U_1^2(t)\tau(\tau-2)(4-\tau^2)}{32(1+\beta)(1+3\beta)} \le 0,$$

$$c_3(t,\tau) = \frac{U_1^3(t)\tau^2(4-\tau^2)}{64(1+\beta)^2(1+2\beta)} + \frac{U_1(t)U_2(t)\tau^2(4-\tau^2)}{16(1+\beta)(1+3\beta)} \ge 0,$$

$$c_4(t,\tau) = \frac{U_1(t)\left|(1+\beta)^3 U_3(t) - (1+3\beta)U_1^3(t)\right|}{16(1+\beta)^4(1+3\beta)}\tau^4 + \frac{U_1^2(t)\tau(4-\tau^2)}{8(1+\beta)(1+3\beta)} \ge 0, t \in \left(\frac{1}{2}, 1\right), \tau \in [0,2].$$

Now, we need to maximize the function $F(\xi,\eta)$ on the closed square $\Omega = \{(\xi,\eta): \xi,\eta \in [0,1]\}$ for $\tau \in [0,2]$. Since the coefficients of the function $F(\xi,\eta)$ is dependent to variable $\tau$ for fixed value of $t$, we must investigate the maximum of $F(\xi,\eta)$ respect to $\tau$ taking into account these cases $\tau = 0$, $\tau = 2$ and $\tau \in (0,2)$.

Let $\tau = 0$. Then, we write

$$F(\xi,\eta) = c_1(t,0) = \frac{U_1^2(t)}{4(1+2\beta)^2}(\xi+\eta)^2.$$

It is clear that the maximum of the function $F(\xi,\eta)$ occurs at $(\xi,\tau) = (1,1)$, and

$$\max\{F(\xi,\eta) : \xi,\eta \in [0,1]\} = F(1,1) = \frac{U_1^2(t)}{(1+2\beta)^2} \quad (2.22)$$

Now, let $\tau = 2$. In this case, $F(\xi,\eta)$ is a constant function (respect to $\tau$) as follows:

$$F(\xi,\eta) = c_4(t,2) = \frac{U_1(t)\left|(1+\beta)^3 U_3(t) - (1+3\beta)U_1^3(t)\right|}{(1+\beta)^4 (1+3\beta)}. \quad (2.23)$$

In the case $\tau \in (0,2)$, we will examine the maximum of the function $F(\xi,\eta)$ taking into account the sing of $\Lambda(\xi,\eta) = F_{\xi\xi}(\xi,\eta)F_{\eta\eta}(\xi,\eta) - \left[F_{\xi\eta}(\xi,\eta)\right]^2$.

By simple computation, we can easily see that

$$\Lambda(\xi,\eta) = 4c_2(t,\tau)\left[2c_1(t,\tau) + c_2(t,\tau)\right].$$

Since $c_2(t,\tau) < 0$ for all $t \in \left(\frac{1}{2},1\right], \tau \in (0,2)$ and

$$2c_1(t,\tau) + c_2(t,\tau) = \frac{U_1^2(t)(4-\tau^2)(2-\tau)}{32(1+\beta)(1+2\beta)^2(1+3\beta)}\varphi(\tau),$$

where $\varphi(\tau) = 2(1+\beta)(1+3\beta) - \beta^2\tau$ since $2 - \beta^2\tau > 0$, $\varphi(\tau) > 6\beta^2 + 8\beta^2 > 0$; that is $2c_1(t,\tau) + c_2(t,\tau) > 0$ for all $\tau \in (0,2)$ and $\beta \in [0,1]$, we conclude that $\Lambda(\xi,\eta) < 0$ for all $(\xi,\eta) \in \Omega$. Consequently, the function $F(\xi,\eta)$ cannot have a local maximum in $\Omega$. Therefore, we must investigate the maximum of the function $F(\xi,\eta)$ on the boundary of the square $\Omega$.

Let
$$\partial\Omega = \{(0,\eta):\eta \in [0,1]\} \cup \{(\xi,0):\xi \in [0,1]\} \cup \{(1,\eta):\eta \in [0,1]\} \cup \{(\xi,1):\xi \in [0,1]\}.$$

We can easily show that the maximum of the function $F(\xi,\eta)$ on the boundary $\partial\Omega$ of the square $\Omega$ occurs at $(\xi,\eta) = (1,1)$, and

$$\max\{F(\xi,\eta) : (\xi,\eta) \in \partial\Omega\} = F(1,1) =$$
$$4c_1(t,\tau) + 2[c_2(t,\tau) + c_3(t,\tau)] + c_4(t,\tau), t \in \left(\frac{1}{2},1\right], \tau \in (0,2). \quad (2.24)$$

Now, let us define the function $H:(0,2)\to\mathbb{R}$ as follows:

$$H(t,\tau)=4c_1(t,\tau)+2[c_2(t,\tau)+c_3(t,\tau)]+c_4(t,\tau) \qquad (2.25)$$

for fixed value of $t$.
Substituting the value $c_j(t,\tau), j=1,2,3,4$ in the (2.25), we obtain

$$H(t,\tau)=\frac{U_1^2(t)}{(1+2\beta)^2}+\frac{\Delta(\beta,t)\tau^4+4c(\beta,t)\tau^2}{32(1+\beta)^4(1+2\beta)^2(1+3\beta)},$$

where

$$\Delta(\beta,t)=2U_1(t)\left|(1+\beta)^3 U_3(t)-(1+3\beta)U_1^3(t)\right|(1+2\beta)^2-$$
$$U_1(t)\left[U_1^2(t)-4U_2(t)(1+\beta)(1+2\beta)\right](1+\beta)^2(1+2\beta)(1+3\beta)-2U_1^2(t)\beta^2(1+\beta)^3,$$
$$c(\beta,t)=U_1(t)\left[U_1^2(t)(1+3\beta)+4U_2(t)(1+\beta)(1+2\beta)\right](1+\beta)^2(1+2\beta)-$$
$$2U_1^2(t)\left[(1+\beta)^2+(2+\beta)\beta\right](1+\beta)^3.$$

Now, we must investigate the maximum (respect to $\tau$) of the function $H(t,\tau)$ in the interval $(0,2)$ for fixed value of $t$.

By simple computation, we can easily show that

$$H'(t,\tau)=\frac{\Delta(\beta,t)\tau^2+2c(\beta,t)}{8(1+\beta)^4(1+2\beta)^2(1+3\beta)}\tau.$$

We will examine the sign of the function $H'(t,\tau)$ depending on the different cases of the signs of $\Delta(t,\tau)$ and $c(t,\tau)$ as follows.

(ı) Let $\Delta(\beta,t)\geq 0$ and $c(\beta,t)\geq 0$, then $H'(t,\tau)\geq 0$, so $H(t,\tau)$ is an increasing function. Therefore,

$$\max\{H(t,\tau):\tau\in(0,2)\}=H(t,2-)=\frac{U_1(t)\left|(1+\beta)^3 U_3(t)-(1+3\beta)U_1^3(t)\right|}{(1+\beta)^4(1+3\beta)}. \qquad (2.26)$$

That is,

$$\max\{\max\{F(\xi,\eta):\xi,\eta\in[0,1]\}:\tau\in(0,2)\}=H(t,2-).$$

(ıı) Let $\Delta(\beta,t)>0$ and $c(\beta,t)<0$, then $\tau_0=\sqrt{\dfrac{-2c(\beta,t)}{\Delta(\beta,t)}}$ is a critical point of the function $H(t,\tau)$. We assume that $\tau_0\in(0,2)$. Since $H''(t,\tau_0)>0$, $\tau_0$ is a local minimum point of the function $H(t,\tau)$. That is, the function $H(t,\tau)$ cannot have a local maximum.

(*iii*) Let $\Delta(\beta,t) \leq 0$ and $c(\beta,t) \leq 0$, then $H'(t,\tau) \leq 0$. Thus, $H(t,\tau)$ is an decreasing function on the interval $(0,2)$.

Therefore,

$$\max\{H(t,\tau): \tau \in (0,2)\} = H(t,0+) = 4c_1(t,0) = \frac{U_1^2(t)}{(1+2\beta)^2}. \qquad (2.27)$$

(*iv*) Let $\Delta(\beta,t) < 0$ and $c(\beta,t) > 0$, then $\tau_0$ is a critical point of the function $H(t,\tau)$. We assume that $\tau_0 \in (0,2)$. Since $H''(t,\tau_0) < 0$, $\tau_0$ is a local maximum point of the function $H(t,\tau)$ and maximum value occurs at $\tau = \tau_0$.

Therefore,

$$\max\{H(t,\tau): \tau \in (0,2)\} = H(t,\tau_0), \qquad (2.28)$$

where

$$H(t,\tau_0) = \frac{4t^2}{(1+2\beta)^2} - \frac{c^2(\beta,t)}{4(1+\beta)^4(1+2\beta)^2(1+3\beta)}.$$

Thus, from (2.22)-(2.28), the proof of Theorem 2.1 is completed.

In the special cases from Theorem 2.1, we arrive at the following results.

**Corollary 2.1.** *Let the function $f(z)$ given by (1.1) be in the class $Q_\Sigma(G,0,t) = N_\Sigma(G,t)$, $t \in \left(\frac{1}{2},1\right)$, where the function $G$ is an analytic function given by (1.4). Then,*

$$\left|a_2 a_4 - a_3^2\right| \leq 8t^2.$$

**Corollary 2.2.** *Let the function $f(z)$ given by (1.1) be in the class $Q_\Sigma(G,1,t) = \Re_\Sigma(G,t)$, $t \in \left(\frac{1}{2},1\right)$, where the function $G$ is an analytic function given by (1.4). Then,*

$$\left|a_2 a_4 - a_3^2\right| \leq t^2(1-t^2).$$

From the Corollary 2.2, we can easily arrive at the following result.

**Corollary 2.3.** *Let the function $f(z)$ given by (1.1) be in the class $Q_\Sigma(G,1,\frac{\sqrt{2}}{2}) = \Re_\Sigma(G)$, where the function $G$ is an analytic function given by (1.4). Then,*

$$|a_2 a_4 - a_3^2| \leq \frac{1}{4}.$$

## Acknowledgement

This work was supported by Kafkas University the Scientific Research Projects Coordination Unit. Also, the author is grateful to the anonymous referees for the valuable comments and suggestions.